\documentclass[a4paper,11pt]{article}

\usepackage{amssymb}
\usepackage{amsthm}

\usepackage{mathrsfs}

\newcommand{\van}{{\mathscr L}_{\scriptscriptstyle{0}}}

\newcommand{\dc}{d_{{\scriptscriptstyle
    c}}}
\newcommand{\bu}{{\mathscr B}_{1}}
\newcommand{\bd}{{\mathscr B}_{2}}
\newcommand{\integer}{{\mathbb N}}

\newtheorem{theorem}{Theorem}[section]
\newtheorem{lemma}{Lemma}[section]
\newtheorem{proposition}{Proposition}[section]
\newtheorem{definition}{Definition}[section]
\newtheorem{hypothesis}{Hypothesis}[section]
\newtheorem{remark}{Remark}[section]

\newcommand{\proba}{{\mathbb P}}
\newcommand{\expectation}{{\mathbb E}}
\newcommand{\Oun}{\mathcal{O}(1)} 
\font\gfont=cmmi10 scaled \magstep2
\newcommand{\gchi}{\hbox{\gfont \char31}}
\newcommand{\ka}[1]{\kappa\big(\mathcal{A}_{#1},\mathcal{N}(0,\sigma^{2})\big)}

\renewenvironment{proof}[1]
{\noindent{\bf Proof.}\hspace{0.1cm} #1} {\hfill$\blacksquare$\bigskip\bigskip}

\begin{document}

\begin{center}
{\Large{\bf Statistical Consequences of\\ Devroye Inequality for Processes.\\
Applications to a Class of Non-Uniformly Hyperbolic Dynamical Systems}}
\end{center}

\bigskip

\begin{center}
J.-R. Chazottes$^{a,}$\footnotemark[1]\ , P. Collet$^{a,}$\footnotemark[1]\ 
and B. Schmitt$^{b,}$\footnotemark[1]\ 
\end{center}

\begin{center}
$^a$Centre de Physique Th{\'e}orique,\\
Ecole polytechnique, CNRS UMR 7644\\
F-91128 Palaiseau Cedex, France\\
emails: {\tt jeanrene@cpht.polytechnique.fr}\\
 {\tt collet@cpht.polytechnique.fr}
\end{center}

\begin{center}
$^b$ D{\'e}partement de Math{\'e}matiques \\
Universit{\'e} de Bourgogne\\
Facult{\'e} des Sciences Mirande\\
BP 138, 21004 Dijon Cedex, France\\
email~: {\tt schmittb@u-bourgogne.fr}
\end{center}

\footnotetext[1]{{\bf Acknowledgments}. JRC and PC acknowledge the CIC
for its warm hospitality, Cuernavaca, M\'exico, where part
of this work has been done, as well as the Institut de Math\'ematiques de
Bourgogne in Dijon. BS acknowledges the kind hospitality of CPhT at
Ecole Polytechnique. The authors acknowledge the anonymous referee for
a very careful reading of the paper.}

\begin{abstract}
In this paper, we apply Devroye inequality to 
study various statistical estimators and fluctuations of observables
for processes. Most of these observables are suggested by dynamical systems. 
These applications concern the co-variance function, the integrated periodogram,
the correlation dimension, the kernel density estimator, the speed of
convergence of the empirical measure, the shadowing property and the
almost-sure central limit theorem.
We proved in \cite{CCS} that Devroye inequality holds for
a class of non-uniformly hyperbolic dynamical systems introduced
in \cite{young}.
In the second appendix we prove that, if the decay of correlations
holds with a common rate for all pairs of functions, then it holds
uniformly in the function spaces.
In the last appendix we prove that for the subclass of one-dimensional
systems studied in \cite{young} the density of the absolutely continuous
invariant measure belongs to a Besov space.

\bigskip

\noindent {\bf Keywords}. Integrated periodogram, correlation
dimension, kernel density estimator, empirical measure, shadowing,
almost-sure central limit theorem.

\end{abstract}

\newpage

\tableofcontents

\section{Introduction and set-up} 

Assume one has a finite sample $x_1,\ldots,x_n$ 
of a stationary ergodic process taking values in
$\mathbb R^d$.
If we consider an empirical estimator (or an observable) $K(x_1,\ldots,x_n)$
of some statistical properties of the process, we basically wish to
determine its fluctuations and its convergence properties, as $n$
grows.
In Statistician's terminology, we aim to study the consistency of the estimator
$K(x_1,\ldots,x_n)$ and be able to build confidence intervals.

As we shall see in the sequel with various examples, many interesting estimators 
have a complicated dependence on the sample. In particular they are
not of the form $(u(x_1)+\cdots+u(x_n))/n$, for some function $u$, or
cannot be well approximated by such time-averages for which the
Central Limit Theorem may apply.

The aim of this paper is to apply what we call Devroye inequality \cite{devroye}, see
the definition below, to estimate the variance for a general class of
estimators $K(x_1,\ldots,x_n)$. 
For some of them we will further require some weak conditions on the
auto-covariance function for functionals of the process. 

Our applications concern the empirical auto-covariance function, the
integrated periodogram, the correlation dimension, the kernel
density estimation of the density of the invariant measure, shadowing
properties, the speed of convergence of the empirical measure toward
the invariant measure, and the almost-sure
central limit theorem. Some of these estimators were
studied in \cite{cms} in the context of piece-wise expanding maps on the
interval for which a stronger inequality than Devroye inequality holds.

We shall formulate the results as much as possible in an abstract
setting in order to see more clearly what is needed to prove them.
As we showed in \cite{CCS}, a class of non-uniformly hyperbolic 
dynamical systems introduced by L.-S. Young
\cite{young} fits this framework.

Let $(\Omega,\mathfrak{B},\proba)$ be a probability space and $(X_k)$
be a stationary ergodic sequence of random variables assuming values
in $\mathbb{R}^d$.

We will denote the expectation with respect to $\proba$ by
$\expectation$, and by $\mu$ the common distribution of
the $X_k$'s. We will assume that the $X_k$'s are 
almost-surely bounded, i.e. there exists a positive constant $A$ such that
\begin{equation}\label{A}
\Vert X_k\Vert \leq A \quad\proba-\textup{almost surely}\,. 
\end{equation}

Let $K$ be a real-valued function on $(\mathbb{R}^{d})^n$. We will say that
$K$ is separately $\eta$-H{\"o}lder in all its variables, if for any
$1\le i\le n$, the following quantities are finite
\begin{equation}\label{Holdercoef}
L_j=L_j(K):=
\sup_{x_1,x_2,\ldots,x_{j-1},x_j,x_{j+1},\ldots,x_n}\,
\sup_{\tilde{x}_j \neq x_j}
\end{equation}
$$
\frac{
|K(x_1,\ldots,x_{j-1},x_{j},x_{j+1},\ldots,x_n)-K(x_1,\ldots,x_{j-1},\tilde{x}_{j},x_{j+1
},\ldots,x_n)|
}{\Vert x_j - \tilde{x}_j \Vert^\eta}\;\cdot
$$

We now define what we mean by saying that the  process $(X_k)$
satisfies Devroye inequality.

\begin{definition}[Devroye inequality for the variance]
We will say that the process $(X_k)$ satisfies Devroye inequality if,
for $\eta\in]0,1]$,
there exists a constant $D=D(\eta)>0$ such that for any integer $n\geq 1$ and
for any real-valued separately $\eta$-H{\"o}lder function $K$ on
$(\mathbb{R}^d)^n$, we have
\begin{equation}\label{DEV}
\textup{var}(K)=\expectation\left(\big( K - \expectation(K)\big)^2\right)
\leq
D\ \sum_{j=1}^{n} L_{j}^2\;.
\end{equation}
\end{definition}

For the case of Dynamical Systems, $\Omega$ is the phase space on
which acts a measurable transformation $f$. We assume that an
$f$-invariant ergodic measure $\mu$ is given. One can define a
stochastic process $X_k(x)=f^{k-1}(x)$ where $x$ is randomly chosen
according to $\mu$. We are interested in observables of the form
$K(X_1,\ldots,X_n)(x)=K(x,f(x),\ldots,f^{n-1}(x))$.

One can ask whether there are processes satisfying
Devroye inequality.
Indeed, a large class of dynamical systems satisfy Devroye
inequality, as we proved in \cite{CCS}.  
Let us recall that this class contains families of
piece-wise hyperbolic maps, like the Lozi maps; scattering billiards,
like the planar periodic Lorentz gas; quadratic and H{\'e}non maps
(for parameter sets with positive Lebesgue measure). 
Let us also briefly mention that such dynamical systems admit an
SRB-measure, enjoy exponential decay of correlations and a central
limit theorem for H\"older continuous observables. 
Notice that in the sequel we will only need very slow decays of
correlations, e.g., $C(\ell)\sim 1/\sqrt{\ell}$ for the integrated
periodogram or absolute summability for, e.g., the almost-sure central
limit theorem.

\section{Covariance function}

Recall that the auto-covariance of a real-valued, square-integrable,
function $u$ on ${\mathbb R}^d$ is defined by
\begin{equation}\label{DC}
C(n)=C_u(n):=\expectation (u(X_1) u(X_n)) - (\expectation(u(X_1)))^2\,.
\end{equation}

An empirical estimator of the auto-covariance is given by
$$
\hat{C}_{k}(n)=\frac{1}{k}\sum_{j=1}^k u(X_j) u(X_{j+n}) -
\left(\frac{1}{k}\sum_{j=1}^k u(X_j)\right)^2
\;.
$$
It follows at once from Birkhoff's ergodic Theorem that
$$
C(n)=\lim_{k\to\infty}\hat{C}_{k}(n)\quad\proba-\textup{almost surely}\,.
$$

\begin{theorem}
Let $u$ be a real-valued $\eta$-H\"older function on ${\mathbb R}^d$ with
H\"older constant denoted by $L_u$. Then, for all integers $k$, $n$, we have
$$
\expectation\left( \big(\hat{C}_{k}(n)-C(n)\big)^2\right)
\leq
16 D L_u^4 A^{2\eta}\frac{n+k}{k^2} + \frac{D^2 L_u^4}{k^2}\,\cdot
$$
\end{theorem}

\begin{proof}
We have the following identity~:
$$
\expectation\left(
  \big(\hat{C}_{k}(n)-C(n)\big)^2\right)=
$$
$$
\expectation\left(
  \big(\hat{C}_{k}(n)-\expectation(\hat{C}_{k}(n))\big)^2\right)+
  \big(\expectation(\hat{C}_{k}(n))-C(n))\big)^2 =
$$
$$
\textup{var}(\hat{C}_{k}(n)) + \left(\textup{var}\Big(\frac{1}{k}\sum_{j=1}^k u(X_j)\Big)\right)^2\,.
$$
The first term is estimated using Devroye inequality
(\ref{DEV}) and assuming, without loss of generality, that $u(0)=0$.
We obtain the upper-bound of independent interest
$$
\textup{var}(\hat{C}_{k}(n))\leq 16 D L_u^4 A^{2\eta}\frac{n+k}{k^2}\,\cdot
$$
The second term is easily estimated using again Devroye inequality. 
This leads immediately to the above estimate.
\end{proof}

\noindent {\bf Remark}. For the study of $U$-statistics of functionals of $\alpha$- and
$\beta$-mixing process we refer the interested reader to \cite{B} and
references therein.

\section{Integrated periodogram}

We recall (see \cite{BD}) that if $u$ is a real-valued function 
the raw periodogram (of order $n$) of the process $(u(X_k))$ is the function
\begin{equation}
I_n(\omega)= \frac1{n}\left| \sum_{j=1}^n e^{-ij\omega}\ \left(u(X_j)-
  \expectation(u(X_1))\right)\right|^2\
\end{equation}
where $\omega\in[0,2\pi]$.
The spectral distribution function of order $n$ (integral of the raw
periodogram of order $n$) is given by
\begin{equation}\label{SDF}
J_n(\omega)= \int_0^\omega I_n(s)\ ds\, .
\end{equation}

From a practical point of view, it is worth defining the empirical 
spectral distribution function of order $n$ as follows:
$$
\tilde{J}_n(\omega)= \int_0^\omega \frac{1}{n}\left| \sum_{j=1}^{n} e^{-ijs}
\left( u(X_j)-\frac1{n}\sum_{\ell=1}^{n} u(X_\ell)\right)\right|^2
\ ds\, .
$$

In this section we will make the following assumption.
\begin{hypothesis}\label{pluie}
The function $u$ is $\eta$-H\"older continuous and its 
auto-covariance function $C(\ell)=C_u(\ell)$ 
satisfies
$$
\sum_{\ell=1}^{\infty} \frac{|C(\ell)|}{\ell}<\infty
$$
(where $C(\ell)$ is defined at (\ref{DC})).
\end{hypothesis}

Let $\hat{C}(\omega)$ be the Fourier cosine transform of the
auto-covariance function, namely
$$
\hat{C}(\omega)=\sum_{k=0}^\infty \cos(\omega k)\
C(k+1) \,.
$$
We will denote by $J(\omega)$ the integral of the following quantity
\begin{equation}\label{IFTAF}
J(\omega)= \int_0^\omega (2\hat{C}(s) - C(1))\ ds=
C(1)\ \omega + 2\sum_{k=1}^\infty \frac{\sin(\omega k)}{k}\
C(k+1) \,.
\end{equation}

We will use the following convenient quantity:
$$
\Delta_n := \frac{2}{n} \sum_{k=1}^{n-1} |C(k+1)| + 2 \sum_{k=n}^\infty \frac{|C(k+1)|}{k}\ \cdot
$$

Observe that $J(\omega+2\pi)= J(\omega) + 2\pi C(1)=J(\omega) +
J(2\pi)$. In order to estimate $J$, it is therefore enough to
restrict to the interval $[0,2\pi]$.

\begin{theorem}\label{Jtheorem}
There exists a positive constant $\Gamma$ 
such that for any function $u$ satisfying Hypothesis \ref{pluie}, and any $n\geq 1$, we have:
$$
\expectation \left(\Big( \sup_{\omega\in [0,2\pi]} \big|
\tilde{J}_n(\omega)-J(\omega) \big| \Big)^2 \right)
$$
$$
\leq 
\Gamma\ \inf_{N\geq 1} \left\{ N
\left[\frac{C(1)^2 + D A^{2\eta}\ L_u^4 (1+\log n)^2}{n} +\Delta_n^2\right]
+ 
\left[\frac{C(1)}{N} + \Delta_N \right]^2 \right\}\ .
$$
\end{theorem}

\noindent {\bf Remark}. If $\Delta_n \leq const/n$, then 
$$
\expectation \left(\Big( \sup_{\omega\in [0,2\pi]} \big|
\tilde{J}_n(\omega)-J(\omega) \big| \Big)^2 \right)
\leq 
\Oun \frac{(1+\log n)^{4/3}}{n^{2/3}}\ \cdot
$$
In particular, if the auto-covariance is absolutely summable, then
$\Delta_n \leq const/n$.

For convergence results in distribution sense of the raw periodogram
for a class of maps on the interval we refer to \cite{LL}.

This theorem is the consequence of two propositions.

\begin{proposition}\label{propJn-J}
For any function $u$ satisfying Hypothesis \ref{pluie}, and any $n\geq 1$, we have:
$$
\expectation \left(\Big( \sup_{\omega\in [0,2\pi]} \big|
  J_n(\omega)-J(\omega) \big| \Big)^2 \right)
\leq 
$$
$$
\inf_{N>1}\left\{
2(N+1) \left(
\frac{(4\pi + 1+\log n)^2 L_u^4 D}{n} \ A^{2\eta}  +
\Delta_n^2 \right) + 
8\pi^2 \left( \frac{C(1)}{N} + \Delta_N\right)^2\right\}\ .
$$
\end{proposition}

\begin{proof}
Let
\begin{equation}\label{Qn}
Q_n=\sup_{\omega\in [0,2\pi]} \big| J_n(\omega)-J(\omega) \big|\,.
\end{equation}
Let $N$ be an integer and define the sequence of numbers $(\omega_p)$
by $\omega_p = 2\pi p / N$ for $p=0,\ldots,N$.
It follows at once from the monotonicity of $J$ and $J_n$ (since they
are integrals of non-negative functions) that
$$
Q_n \leq \max\left( \sup_{0\leq p \leq N-1} | J_n(\omega_{p+1})-J(\omega_p)|,
\sup_{0\leq p \leq N-1} | J_n(\omega_{p})-J(\omega_{p+1})|    \right)\,.
$$
We now have
$$
Q_n \leq  \sup_{0\leq p \leq N} | J_n(\omega_{p})-J(\omega_p)|+
\sup_{0\leq p \leq N-1} | J(\omega_{p})-J(\omega_{p+1})|\ .
$$
Now using definition (\ref{IFTAF}), we get after an easy computation that
for all $p=0,\ldots,N-1$
\begin{equation}\label{Jomegai}
| J(\omega_{p})-J(\omega_{p+1})|\leq 2\pi\left( \frac{C(1)}{N} + \Delta_{N}\right)\ . 
\end{equation}
It follows that
\begin{equation}\label{QnQnbar}
Q_n\leq \overline{Q}_n + 2\pi\left( \frac{C(1)}{N} + \Delta_{N}\right)
\end{equation}
where
$$
\overline{Q}_n=\sup_{0\leq p \leq N} | J_n(\omega_{p})-J(\omega_p)|\,. 
$$
We obviously have
\begin{equation}\label{obviousformula}
\expectation(\overline{Q}^2_n )\leq   \sum_{p=0}^{N}
\expectation\big(\left( J(\omega_p) - J_n(\omega_p)  \right)^2 \big)\,.
\end{equation}
We now estimate each term $\expectation\big(\left(J(\omega_p) - J_n(\omega_p)  \right)^2 \big)$.
Observe that for any $\omega$ we have
$$
\expectation\big(\left(J(\omega) - J_n(\omega)  \right)^2 \big)=
\expectation\big(\left(J_n(\omega)-\expectation(J_n(\omega))  \right)^2 \big)+
\left(\expectation(J_n(\omega) - J(\omega))  \right)^2\ . 
$$
We have also from the definition of $J_n$
$$
J_n(\omega)=\frac{\omega}{n}\sum_{j=1}^n (u(X_j)-\expectation(u(X_1))^2 
$$
$$
+\ \frac{i}{n}\sum_{j\neq \ell}^n \frac{e^{-i(j-\ell)\omega}-1}{j-\ell}\
(u(X_j)-\expectation(u(X_1))(u(X_\ell)-\expectation(u(X_1))=
$$
$$
\frac{\omega}{n}\sum_{j=1}^n (u(X_j)-\expectation(u(X_1))^2 
$$
\begin{equation}\label{formulaforJn}
+\ \frac{1}{n}\sum_{j\neq \ell}^{n}
\frac{\sin((j-\ell)\omega)}{j-\ell}\ 
(u(X_j)-\expectation(u(X_1))(u(X_\ell)-\expectation(u(X_1))\,.
\end{equation}
Using this formula and (\ref{IFTAF}), an easy computation leads to
\begin{equation}\label{formulaforEJn}
\left(\expectation(J_n(\omega)) - J(\omega))  \right)^2 \leq \Delta_n^2\ .
\end{equation}
We now apply Devroye inequality to $J_n(\omega)$ in the form
(\ref{formulaforJn}) and get
\begin{equation}\label{otherformulaforEJn-J}
\expectation\big(\left(J_n(\omega)-\expectation(J_n(\omega))  \right)^2 \big)
\leq
\frac{(4\pi + 1+\log n)^2 L_u^4}{n}\ A^{2\eta} D \ .
\end{equation}
Using (\ref{formulaforEJn}), (\ref{otherformulaforEJn-J}) and (\ref{obviousformula}),
it follows that
\begin{equation}\label{EQnbar}
\expectation(\overline{Q}^2_n) \leq (N+1) \left( 
\frac{(4\pi + 1+\log n)^2 L_u^4}{n} A^{2\eta} D  +
\Delta_n^2 \right)\ .
\end{equation}

This completes the proof.
\end{proof}

\begin{proposition}\label{propJn-Jntilde}
There exists a positive constant $S$ 
such that for any function $u$ satisfying Hypothesis \ref{pluie}, and any $n\geq 1$, we have
$$
\expectation \left(\Big( \sup_{\omega\in [0,2\pi]} \big|
  J_n(\omega)-\tilde{J}_n(\omega) \big| \Big)^2 \right)
$$
$$
\leq 
S\ \inf_{N\geq 1}
\left\{
(N+1) \left[\frac{C(1)^2 + D A^{2\eta} L_u^4}{n} +\Delta_n^2\right]
+ 
\left[\frac{C(1)}{N} + \Delta_N \right]^2 \right\}\ .
$$
\end{proposition}

The proof is rather similar to the previous one.

\bigskip

\begin{proof}
Let
\begin{equation}\label{Rn}
R_n=\sup_{\omega\in [0,2\pi]} \big| \tilde{J}_n(\omega)-J_n(\omega) \big|\,.
\end{equation}
Let $N$ be an integer and define as before the sequence of numbers $(\omega_p)$
by $\omega_p = 2\pi p / N$ for $p=0,\ldots,N$.
It follows at once from the monotonicity of $J_n$ and $\tilde{J}_n$ that
$$
R_n \leq \max\left( \sup_{0\leq p \leq N-1} | J_n(\omega_{p+1})-\tilde{J}_n(\omega_p)|,
\sup_{0\leq p \leq N-1} | J_n(\omega_{p})-\tilde{J}_n(\omega_{p+1})|    \right)\,.
$$
We now have
$$
R_n \leq  \overline{R}_n +
\sup_{0\leq p \leq N-1} | J_n(\omega_{p})-J_n(\omega_{p+1})|
$$
where
\begin{equation}\label{Rnbar}
\overline{R}_n=\sup_{0\leq p \leq N} | \tilde{J}_n(\omega_{p})-J_n(\omega_p)|\,.
\end{equation}
Now we have the estimate
$$
\expectation\left( \sup_{0\leq p \leq N-1} (J_n(\omega_p)-J_n(\omega_{p+1}))^2\right)\leq
$$
$$
6\ \expectation\left( \sup_{0\leq p \leq N} (J_n(\omega_p)-J(\omega_{p}))^2 \right)
+
3\ \sup_{0\leq p \leq N-1} (J(\omega_p)-J(\omega_{p+1}))^2 \,.
$$
Using Proposition \ref{propJn-J} to estimate the first term and (\ref{Jomegai})
for the second one, we obtain
$$
\expectation\left( \sup_{0\leq p \leq N-1} (J_n(\omega_p)-J_n(\omega_{p+1}))^2\right)\leq
$$
\begin{equation}\label{EsupJnomegai}
12(N+1) \left(
\frac{4\pi^2 L_u^4 D}{n} A^{2\eta}  +
\Delta_n^2 \right) + 
60\pi^2 \left( \frac{C(1)}{N} + \Delta_N\right)^2\,.
\end{equation}
We obviously have
\begin{equation}\label{obviousformulabis}
\expectation(\overline{R}^2_n )\leq   \sum_{p=0}^{N-1}
\expectation\left( (J_n(\omega_p) - \tilde{J}_n(\omega_p))^2 \right)\,.
\end{equation}
We now have to estimate each term $\expectation\left((J_n(\omega_p) -
\tilde{J}_n(\omega_p))^2  \right)$.
Observe that for any $\omega$
$$
\expectation\left((J_n(\omega) - \tilde{J}_n(\omega))^2\right)=
\textup{var}\left( J_n(\omega) - \tilde{J}_n(\omega)\right) +
\left( \expectation( J_n(\omega) - \tilde{J}_n(\omega))  \right)^2\ .
$$
Let $S_n:=\sum_{j=1}^n u(X_j)$.
A simple computation yields
$$
J_n(\omega)-\tilde{J}_n(\omega)=\omega\ \left(\frac{S_n}{n}-\expectation(u(X_1))\right)^2
$$
\begin{equation}\label{formulaforJn-Jntilde}
+
\ \frac{1}{n}\left(\frac{S_n}{n}-\expectation(u(X_1))\right)
\sum_{j\neq \ell}^n \frac{\sin((j-\ell)\omega)}{j-\ell}
\left(2u(X_\ell)-\expectation(u(X_1))-\frac{S_n}{n}\right)\ .
\end{equation}
An easy computation leads to
$$
\expectation\left(J_n(\omega)-\tilde{J}_n(\omega)\right)=
\left( \omega -\frac1{n} \sum_{j\neq \ell}^n \frac{\sin((j-\ell)\omega)}{j-\ell}\right)
\expectation\left(\Big(\frac{S_n}{n}-\expectation(u(X_1))\Big)^2\right)
$$
$$
+\frac{2}{n^2} \sum_{r=1}^n \sum_{\ell=1}^n
\expectation\bigg( 
\big( u(X_r)-\expectation(u(X_1)) \big)
\big( u(X_\ell)-\expectation(u(X_1)) \big) \bigg)\
 \sum_{j\neq \ell}^n \frac{\sin((j-\ell)\omega)}{j-\ell}\,\cdot
$$
An easy computation using Lemma \ref{sin} shows that there is a constant
$c_1>0$ such that for all integer $n$
$$
\sup_{\omega\in[0,2\pi]}
\left| \left( \omega -\frac1{n} \sum_{j\neq \ell}^n \frac{\sin((j-\ell)\omega)}{j-\ell}\right)
\expectation\left(\Big(\frac{S_n}{n}-\expectation(u(X_1))\Big)^2\right) \right|
\leq 
$$
$$
c_1 \ \left( \frac{C(1)}{n}+\Delta_n \right)\ .
$$
Similarly, there exists a constant $c_2>0$ such that
$$
\sup_{\omega\in[0,2\pi]}
\left|\frac{2}{n^2} \sum_{r=1}^n \sum_{\ell=1}^{n} C(|\ell-r|+1)
\sum_{j\neq \ell}^n \frac{\sin((j-\ell)\omega)}{j-\ell}\right|
\leq c_2 \ \left( \frac{C(1)}{n}+\Delta_n \right)\ .
$$
Combining these two estimates, one gets
\begin{equation}\label{formulaforEJn-Jntilde}
\sup_{\omega\in[0,2\pi]}
\left( \expectation( J_n(\omega) - \tilde{J}_n(\omega))  \right)^2
\leq
c_3\ \left( \frac{C(1)}{n}+\Delta_n \right)^2
\end{equation}
where $c_3>0$ is a constant (independent of $n$).
We now apply Devroye inequality to $J_n(\omega)-\tilde{J}_n(\omega)$ using (\ref{formulaforJn-Jntilde})
and Lemma \ref{sin}. We easily obtain the estimate
\begin{equation}\label{varJn-Jntilde}
\sup_{\omega\in[0,2\pi]} \ \textup{var}\left(J_n(\omega)-\tilde{J}_n(\omega) \right)
\leq
\frac{c_4 D A^{2\eta} L_u^4}{n}\, \cdot
\end{equation}
It follows that 
$$
\expectation(\overline{R}_n^2) \leq 
N\left( c_3\ \left( \frac{C(1)}{n}+\Delta_n \right)^2+ 
\frac{c_4 D A^{2\eta} L_u^4}{n}\right)\ .
$$
The Proposition follows by combining this estimate with
(\ref{EsupJnomegai}).
\end{proof}

\bigskip

Theorem \ref{Jtheorem} is proved by combining Propositions \ref{propJn-J} and \ref{propJn-Jntilde}.

\section{Correlation dimension}

We recall that the correlation dimension $\dc=\dc(\mu)$
of the measure $\mu$ (recall
that $\mu$ is the common distribution of the $X_k$'s)
is defined by
$$
\lim_{\epsilon\downarrow 0} \frac{\log \int \mu(B(x',\epsilon))\ d\mu(x')}{\log\epsilon^{-1}}
$$
provided the limit exists (where $B(x',\epsilon)$ is the ball of
centre $x'$ and radius $\epsilon$).
In practice one determines for large $n$ the power-law behaviour in
$\epsilon$ of $K_{n,\epsilon}^\vartheta(x,f(x),\ldots,f^{n-1}(x))$ where
$$
K_{n,\epsilon}^\vartheta(x_1,\ldots,x_n)=
\frac1{n^2} \sum_{i\neq j} \vartheta (\epsilon - d(x_i,x_j))
$$
and $\vartheta$ is the Heaviside function (i.e., the characteristic
function of ${\mathbb R}^+$). 
It is known that (see e.g. \cite{MS})
$$
\lim_{n\to\infty}K_{n,\epsilon}^\vartheta(x,f(x),\ldots,f^{n-1}(x))= \int \mu(B(x',\epsilon))\ d\mu(x')
$$
for $\mu$-almost all $x$ and every continuity point of the non-increasing function
$\epsilon\mapsto \int \mu(B(y,\epsilon))\ d\mu(y)$.

To proceed we need to replace $K_{n,\epsilon}^\vartheta(x_1,\ldots,x_n)$
by a component-wise Lipschitz function. 
For any real-valued Lipschitz function $\phi$, define the sequence of component-wise Lipschitz
functions
\begin{equation}\label{regu}
K_{n,\epsilon}^\phi(x_1,\ldots,x_n):=\frac1{n^2} \sum_{i\neq j}
\phi\left(1- \frac{d(x_i,x_j)}{\epsilon}\right)\,\cdot 
\end{equation}

\begin{theorem}
For any real-valued Lipschitz function $\phi$, for any $0<\eta\leq 1$,
there exists a constant $C=C(\eta)>0$ such that for any
$\epsilon>0$ and any integer $n$, we have
\begin{equation}
\textup{var}(K_{n,\epsilon}^\phi) \leq \frac{C}{\epsilon^{2\eta} n}\;\cdot
\end{equation}
\end{theorem}

The proof is a direct application of Devroye inequality (\ref{DEV}).

Several functions $\phi$ are used in the literature. A simple one is given by
$$
\phi_0(y)=\left\{
\begin{array}{l}
0\quad \textup{for}\; y<-\frac1{2}\\
\frac1{2}+y \quad\textup{for}\; -\frac1{2}\leq y\leq\frac1{2}\\
1\quad\textup{for}\;y>\frac1{2}\,\cdot
\end{array}
\right.
$$
One verifies easily that for all $y\in\mathbb{R}$
\begin{equation}\label{regu-example}
\vartheta(1-2y)\leq \phi_0(1-y)\leq \vartheta(1-y/2)\,.
\end{equation}

This implies immediately 
\begin{equation}\label{thetavarphi0}
K_{n,\epsilon/2}^\vartheta(x_1,\ldots,x_n)\leq K_{n,\epsilon}^{\phi_0}(x_1,\ldots,x_n)\leq
K_{n,2\epsilon}^\vartheta(x_1,\ldots,x_n)
\end{equation}
for all $x_1,\ldots,x_n$, $\epsilon>0$ and $n\geq 1$.
It follows that, when $\dc>0$, we have 
$$
K_{n,\epsilon}^\vartheta(x,f(x),\ldots,f^{n-1}(x))\approx \epsilon^{\dc}\quad\textup{as}\;\epsilon\to 0
$$
is equivalent to
$$
K_{n,\epsilon}^{\phi_0}(x,f(x),\ldots,f^{n-1}(x))\approx \epsilon^{\dc}\quad\textup{as}\;\epsilon\to 0\,.
$$
Requiring that the typical value is smaller than the size of fluctuations (standard deviations) leads
to $\epsilon^{\dc} \gtrsim 1/(\epsilon^{\eta}\sqrt{n})$. In other words
$$
n\gtrsim \epsilon^{-2(\dc+\eta)}\,.
$$
In some iid cases, the optimal estimate has been obtained in \cite{optimal}.

\section{Empirical measure}\label{dijon}

We recall that the empirical measure of a sample $X_1,\ldots,X_n$ is a
random measure on $\mathbb{R}^d$ defined by
$$
\mathcal{E}_n= \frac{1}{n} \sum_{j=1}^{n} \delta_{X_j}
$$
where $\delta$ denotes the Dirac measure. 
We recall that from Birkhoff's ergodic theorem, 
almost-surely this sequence of random measures weakly converges to the common
distribution $\mu$ of the $X_k$'s. It is natural to ask for 
the speed of this convergence. This of course depends on the distance
chosen on the set of probability measures. We will consider the
Kantorovich distance defined for two probability measures $\mu_1$ and
$\mu_2$ on ${\mathbb R}^d$ by
\begin{equation}\label{marron}
\kappa(\mu_1,\mu_2)=\sup_{g\in{\mathscr L}}\int g(x)\;d\big(\mu_{1}-\mu_{2}\big)(x)
\end{equation}
where ${\mathscr L}$ denotes the set of real-valued Lipschitz functions
on ${\mathbb R}^d$ with Lipschitz constant at most one.

We now state the theorems of this section. 

\begin{theorem}\label{barque}
By Devroye inequality (\ref{DEV}) we have, for all $n\geq 1$,
$$
\textup{var}(\kappa(\mathcal{E}_n,\mu))\leq \frac{D(1)}{n}\,\cdot
$$
\end{theorem}

The proof follows at once from Devroye inequality (\ref{DEV}) using
the following separately Lipschitz function of $n$ variables
$$
K(x_1,\ldots,x_n)=\sup_{g\in{\mathscr L}}\left[
  \frac{1}{n}\sum_{j=1}^n g(x_j) - \expectation(g)
\right] \,.
$$

To get a probability estimate based on this result one needs to give
an upper-bound for $\expectation(\kappa(\mathcal{E}_n,\mu))$. 
The bound we are so far able to obtain in dimension larger $1$ is too
pessimistic. We explain below how to obtain a more satisfactory bound
in dimension $1$. We will require the following property for the auto-covariance.
We will denote by $\Vert u \Vert_\eta$ the $\eta$-H\"older constant of
$u$ (which is bounded by $\Oun L_1(u)$).

\begin{hypothesis}\label{nuage}
For any $\eta\in]0,1]$ there is a constant $C_{\eta}>0$ such that the 
auto-covariance $C_u(\ell)$ of any $\eta$-H\"older continuous
function $u$ satisfies
$$
\sum_{\ell=1}^{\infty} |C_u(\ell)| \leq C_\eta \ \Vert u\Vert_\eta^2\,.
$$
\end{hypothesis}

This leads to the following theorem.

\begin{theorem}\label{gondole}
Assume that the process $(X_k)$ takes values in $\mathbb{R}$ 
and that the auto-covariance of $\eta$-H\"older continuous
functions satisfies Hypothesis \ref{nuage}.
Then, for any $\eta\in]0,1]$, there exists a positive constant $a(\eta)$ such that
for all $t>0$ and $n\geq 1$, we have
$$
\proba\left(
\kappa(\mathcal{E}_n,\mu)> t+ \frac{a(\eta)}{n^{1/(2(1+\eta))}}\right)\leq
\frac{D(1)}{n t^2} \,\cdot
$$
\end{theorem}

\bigskip

\noindent {\bf Remark}. If $a(\eta)$ behaves like $1/\eta$ as $\eta$
tends to zero, then one can optimize by taking $\eta=1/\log n$.

\bigskip

\begin{proof}
The theorem of Dall'Aglio \cite{dallaglio} states that
$$
\kappa(\mu_1,\mu_2)=\int_{{\mathbb R}} \left|
F_{\mu_1}(t)-F_{\mu_2}(t)\right| \ dt
$$
where $F_{\mu}(t)$ is the distribution function of $\mu$.

We wish to estimate the Kantorovich distance between the empirical 
measure $\mathcal{E}_n$ and $\mu$ (the common distribution of the
$X_k$'s).
In this case we have
$$
\kappa(\mathcal{E}_n,\mu)= \int_{-A}^{A} dt \left|
\frac{1}{n}\sum_{k=0}^{n-1} \vartheta(t-X_k)- F_{\mu}(t)\right|
$$
since we assumed from the very beginning that $\Vert X_k\Vert
\leq A$ $\proba$-almost-surely, and $\vartheta$ denotes the Heaviside function.

In order to use the decay of correlations, we replace the Heaviside
function by a H\"older continuous function $g_\delta$ parametrised by
a positive $\delta$ and defined by
$$
g_\delta(s)=\left\{
\begin{array}{l}
0\quad\quad\quad\quad   \textup{if}\quad s<-\delta \\
1 + s/\delta \quad\; \textup{if}\quad -\delta\leq s\leq 0\\
1 \quad\quad\quad\quad  \textup{if}\quad s>0\,.
\end{array}\right.
$$
We immediately obtain 
\begin{equation}\label{poire}
\kappa(\mathcal{E}_n,\mu)\leq \delta  +
\int_{-A}^{A} dt \left|
\frac{1}{n}\sum_{k=0}^{n-1} g_\delta(t-X_k)- F_{\mu}(t)\right|\,.
\end{equation}
We have
$$
\expectation(\kappa(\mathcal{E}_n,\mu))\leq
\delta + 
\expectation\left(
\int_{-A}^{A} dt \left|
\frac{1}{n}\sum_{k=0}^{n-1} g_\delta(t-X_k)- \expectation(g_\delta(t-X_1))\right|\right)
\; +
$$
$$
\int_{-A}^{A} dt \ \expectation\left| g_\delta(t-X_1)-\vartheta(t-X_1)\right|\leq
$$
$$
2\delta + \expectation\left(
\int_{-A}^{A} dt \left|
\frac{1}{n}\sum_{k=0}^{n-1} g_\delta(t-X_k)- \expectation(g_\delta(t-X_1))\right|\right)\,.
$$
Using Cauchy-Schwarz inequality as in \cite{cms}, one is led to use
the decay of auto-covariance of the functions $g_\delta(t-\cdot)$.
Using Hypothesis \ref{nuage} we get
$$
\expectation(\kappa(\mathcal{E}_n,\mu))\leq 2\delta + \frac{\Oun}{\delta^{\eta}\sqrt{n}}\,\cdot
$$
Using Chebychev inequality, the above estimate with
$\delta=n^{-1/2(1+\eta)}$ and Theorem \ref{barque} we get
the theorem.
\end{proof}

For the application to dynamical systems satisfying Devroye inequality
(see \cite{CCS}), we need moreover to verify Hypothesis \ref{nuage}.
It is often proved, see e.g. \cite{young}, that the auto-covariance
of observables belonging to a Banach space have a common upper bound 
for their rate of decay. It turns out that this implies a uniform
rate of decay for all functions of norm less than or equal to one,
this is the content of Theorem \ref{collet} proved in Appendix \ref{UDC}.
So, if this decay is summable then Theorem \ref{gondole} holds.
For the systems studied in \cite{young}, Hypothesis \ref{nuage}
can be deduced using the estimates provided by approximations \#1
and \#2 and point 4.2. appearing in that paper.

\section{Kernel density estimation for 1D maps}

In this section we assume that $d=1$, namely that the process $(X_k)$ takes values in a bounded
interval of $\mathbb{R}$. Moreover we assume that the common distribution $\mu$ of the $X_k$'s
is absolutely continuous (with respect to Lebesgue measure) and denote
by $\Phi$ its density.
We consider the random empirical densities $(h_n)$ defined by
$$
h_n(X_1,\ldots,X_n;s)= \frac{1}{n\alpha_n} \sum_{j=1}^{n} \psi((s-X_j)/\alpha_n)
$$
where $\alpha_n$ is a positive sequence converging to $0$ and such
that $n\alpha_n$ converges to $+\infty$, and 
$\psi$ (the kernel) is a bounded, non-negative,
Lipschitz continuous function with compact support whose integral
equals $1$. We are interested in the $L^1$ convergence of these
empirical densities to the density $\Phi$ of the common distribution $\mu$ of the $X_k$'s.

\begin{theorem}\label{graine}
Assume that the probability density $\Phi$ satisfies
\begin{equation}\label{moutarde}
\int |\Phi(s) - \Phi(s-y)| \ ds \leq C\ |y|^\tau 
\end{equation}
for some $C>0$, $\tau>0$ and any $y\in\mathbb{R}$. Suppose also that
Hypothesis \ref{nuage} holds.
Then, for any $\eta\in]0,1]$, for any $\psi$ as above, there exists a constant $C'=C'(\eta,\psi)>0$
such that for any integer $n$ and for any $t>C'(\alpha_n^\tau + 1/(\sqrt{n}\alpha_n^{1+\eta}))$, we have
$$
\proba\left( \int |h_n(X_1,\ldots,X_n;s)-\Phi(s)|\ ds > t \right) \leq
\frac{C'}{t^2 n \alpha_n^{2\eta}}\,\cdot 
$$
\end{theorem}

\begin{proof}
We define the functions 
$$
K(x_1,\ldots,x_n)= \int \left| \frac{1}{n\alpha_n} \sum_{j=1}^{n}
  \psi((s-x_j)/\alpha_n) -\Phi(s) \right|\ ds \,. 
$$
It is easy to verify that the H\"older constants of this $\eta$-H\"older continuous function satisfy
$$
\max_{1\leq j\leq n} L_j \leq \frac{\Oun}{n\alpha_n^\eta}\,\cdot
$$
Hence, using Devroye inequality (\ref{DEV}), we immediately obtain
$$
\textup{var}(K)\leq \frac{\Oun}{n\alpha_n^{2\eta}}\,\cdot
$$
The theorem will follow using this and Chebychev inequality provided we have an upper
bound for $\expectation(K)$. To this purpose we will follow the lines
of the proof of Theorem III.2 in \cite{cms} 
with the appropriate modifications.

We first estimate the $L^1$-norm of $\Phi - \expectation(h_n)$. We
obtain, using (\ref{moutarde}) the upper bound 
$$
\int |\Phi(s) - \expectation(h_n)(s)|\ ds \leq \alpha_n^{-1}\int
\psi(y/\alpha_n) \ dy \ 
\int |\Phi(s) - \Phi(s-y)| \ ds \leq
$$
$$
C\ \alpha_n^{-1}\ \int \psi(y/\alpha_n)\ |y|^\tau \ dy \leq
\Oun\  \alpha_n^{\tau}\,.
$$
We now bound from above the integral
$$
\int ds\ \expectation\left( |h_n(X_1,\ldots,X_n;s) - \expectation(h_n)(s)| \right)\,.
$$
By a well-known computation we have
$$
\textup{var}(h_n(X_1,\ldots,X_n;s)) \leq 
$$
$$
\frac{2}{n\alpha_n^2} \sum_{\ell=1}^n 
\expectation\left(
\big(\psi((s-X_1)/\alpha_n)) - \tilde\psi_n(s)\big)
\big(\psi((s-X_\ell)/\alpha_n)) -\tilde\psi_n(s)\big)\right)
$$
where $\tilde\psi_n(s)=\expectation(\psi((s-X_1)/\alpha_n)))$.
Using Cauchy-Schwarz inequality and Hypothesis \ref{nuage}, as in
the proof of Section \ref{dijon}, we get
$$
\int ds\ \expectation\left( |h_n(X_1,\ldots,X_n;s) - \expectation(h_n)(s)| \right) \leq
\frac{\Oun}{\sqrt{n}\ \alpha_n^{1+\eta}}\,\cdot
$$
Summarising we obtain
$$
\expectation(K)\leq \Oun \left(\alpha_n^\tau + \frac{1}{\sqrt{n}\ \alpha_n^{1+\eta}}\right)\,\cdot
$$
The theorem  now follows by Chebychev inequality and Devroye inequality.
\end{proof}

For results on kernel density estimation in the context of piece-wise
expanding maps on the interval, we refer to \cite{clementine} and
references therein.

We will prove in Appendix \ref{bing} that the class of dynamical
systems considered in \cite{young,CCS},
that is the class introduced in \cite{young}, satisfies
(\ref{moutarde}), for 1D systems. 
As explained at the end of the previous section, it also satisfies
Hypothesis \ref{nuage}. Hence the theorem applies.
This class includes quadratic maps for a set of parameter of positive
Lebesgue measure \cite{young}.

\section{Shadowing and mismatch}

For a fixed integer $n$, let $E$ be a measurable subset of $\mathbb{R}^{nd}$. 
For a trajectory $Y_1,\ldots, Y_n$ of length $n$ of
the process $(X_k)$ which is outside $E$, how well can we approximate this trajectory by a trajectory
$(X_1,\ldots,X_n)$ of the process belonging to $E$? 
We first start with a result about the average quality of this
``shadowing''. We will denote by $\mathcal{T}_n$ the set of
trajectories of length $n$ of the process.

\begin{theorem}
For any integer $n$, for any measurable subset $E$ of $\mathbb{R}^{nd}$,
with $\proba((X_1,\ldots, X_n)\in E)>0$, the function
\footnote{The function ${\mathcal Z}_E$ is measurable, see \cite{cms}}
defined by
$$
\mathcal{Z}_E(Y_1,\ldots,Y_n)= \frac{1}{n} 
\inf_{(X_1,\ldots, X_n)\in E\cap \mathcal{T}_n} \sum_{j=1}^{n} \Vert X_j - Y_j \Vert
$$
satisfies for any $t>0$ the inequality
$$
\proba\left( \mathcal{Z}_E \geq \frac{1}{n^{1/3}} \left(t +
    \frac{2^{4/3} D^{1/3}}{\proba((X_1,\ldots, X_n)\in E)}
  \right)\right) 
\leq \frac{D}{n^{1/3} t^2}
$$
where $D>0$ is the constant appearing in (\ref{DEV}).
\end{theorem}

\begin{proof}
We first apply Devroye inequality (\ref{DEV}) to the function
$$
K(x_1,\ldots,x_n)=\frac{1}{n} \inf_{(X_1,\ldots, X_n)\in E\cap
  \mathcal{T}_n} \sum_{j=1}^{n} \Vert X_j - x_j \Vert \,. 
$$
We get
\begin{equation}\label{pouic}
\textup{var}(\mathcal{Z}_E)\leq \frac{D}{n}\,\cdot
\end{equation}
Chebychev inequality yields for any $s>0$
$$
\proba\left( \mathcal{Z}_E \geq \expectation(\mathcal{Z}_E) + \frac{s}{n^{1/3}}\right) \leq
\frac{D}{n^{1/3} s^2}\,\cdot
$$
Proceeding as in \cite{cms} and optimizing over $s$ we obtain
$$
\expectation(\mathcal{Z}_E) \leq \frac{2^{4/3} D^{1/3}}{n^{1/3} \proba((X_1,\ldots, X_n)\in E) }\,\cdot
$$
The theorem follows using again Chebychev inequality. 
\end{proof}

\begin{remark}
There is another way of estimating from above
$\expectation(\mathcal{Z}_E)$. For this observe that
$\mathcal{Z}_E$ vanishes in $E$. Therefore it follows
from (\ref{pouic}) that
$$
\proba((X_1,\ldots, X_n)\in E)\ (\expectation(\mathcal{Z}_E))^2 \leq \frac{D}{n}\,\cdot
$$
Hence 
$$
\expectation((\mathcal{Z}_E)^2) \leq \frac{\sqrt{D}}{\sqrt{n \proba((X_1,\ldots, X_n)\in E)}}\,\cdot
$$
\end{remark}

\bigskip

We now derive a similar result for the number of mismatch at a given precision. 

\begin{theorem}
For any integer $n$, for any measurable subset $E$ of $\mathbb{R}^{nd}$,
with $\proba((X_1,\ldots, X_n)\in E)>0$, and for any $\epsilon>0$, the function defined by
$$
\mathcal{Z}'_{E,\epsilon}(Y_1,\ldots,Y_n)= \frac{1}{n} 
\inf_{(X_1,\ldots, X_n)\in E\cap \mathcal{T}_n} \textup{Card}\{ 1\leq
j\leq n \,:\, \Vert X_j - Y_j \Vert >\epsilon\}  
$$
satisfies for any $t>0$ the following
$$
\proba\left( \mathcal{Z}'_{E,\epsilon} \geq \frac{1}{\epsilon^{2/3} n^{1/3}} 
\left(t + \frac{2^{4/3} D^{1/3}}{\proba((X_1,\ldots, X_n)\in E)} \right)\right)
\leq \frac{D}{\epsilon^{2/3} n^{1/3} t^2}
$$
where $D>0$ is the constant appearing in (\ref{DEV}).
\end{theorem}

\bigskip

The industrious reader can follow the lines of the proof of Theorem
IV.2 in \cite{cms}. Using H\"older estimates 
instead of Lipschitz estimates yields the same formula with
$\epsilon^{2/3}$ replaced with $\epsilon^{2\eta/3}$, for 
any $0<\eta\leq 1$. However the constant $D$ depends on $\eta$ in an implicit way, so it is not clear
how to optimize over $\eta$.

For the case of dynamical systems, given an initial condition $y$ outside
a measurable subset $S$ of the phase space with positive measure, the
questions considered above mean that we look how good is
the shadowing of the orbit of $y$ by an orbit starting from $S$ (in
that case $E=S\times {\mathbb R}^{(n-1)d}$).

\section{Almost-sure central limit theorem}

We say that the process $(u(X_k))$, where $u$ is a real-valued
function, satisfies the Central Limit Theorem if
\begin{equation}\label{clt} 
\lim_{n\to\infty} \proba\left( \frac{\sum_{j=1}^{n} u(X_j) -
    n\expectation(u)}{\sigma\sqrt{n}} \leq t\right) 
= \frac{1}{\sqrt{2\pi}} \int_{-\infty}^{t} e^{-\xi^2/2}\ d\xi
\end{equation}
where
$\sigma^2=\sigma^2(u)$ is assumed to be strictly positive and 
is defined by
\begin{equation}\label{sigma}
\sigma^2= C(1) + 2 \sum_{\ell=2}^{\infty} C(\ell)
\end{equation}
where we assume that the series is finite (see (\ref{DC}) for the
definition of $C(\ell)$).

We will prove an Almost-sure Central Limit Theorem, see e.g. \cite{be}
for a review of this field. Our result is slightly stronger since
it asserts the convergence in the Kantorovich distance
$\kappa$ already used above, see formula (\ref{marron}). We shall use it for
measures on $\mathbb R$ and real-valued Lipschitz functions on
$\mathbb R$.
Note that we can replace $g$ by $g-g(0)$ in (\ref{marron}) since 
$\mu_{1}$ and $\mu_{2}$ are probability measures. In other words there
is no loss of generality in assuming
$$
g\in \van:=\{g\in {\mathscr L}\ | \ g(0)=0\}\,. 
$$
It is convenient to define the sequence of weighted empirical (random)
measures of the normalized partial sum $S_k=u(X_1) + \cdots + u(X_k)$ by 
$$
\mathcal{A}_{n}=\frac{1}{D_n}\sum_{k=1}^{n}\frac{1}{k}
\;\delta_{S_{k}/\sqrt k}
$$
where $D_n=\sum_{k=1}^{n}\frac{1}{k}$. We shall
investigate the convergence of this sequence of weighted empirical
measures to the Gaussian measure in the
Kantorovich metric. 

We now state the result of this section.

\begin{theorem}\label{t1}
Consider the process 
$(u(X_{k}))$ where $u$ is a H\"older continuous function with zero $\mu$
average (recall that $\mu$ is the common law of the $X_k$'s).
Assume that $\sigma^{2}\neq 0$ (see (\ref{sigma})), that the auto-covariance of $(u(X_k))$
is absolutely summable and that (\ref{clt}) holds (central limit theorem).
Then $\proba$-almost surely
\begin{equation}\label{clt1}
\lim_{n\to\infty}\kappa\left(\mathcal{A}_{n},
\mathcal{N}\big(0,\sigma^{2}\big)\right)=0
\end{equation}
where $\mathcal{N}\big(0,\sigma^{2}\big)$ is the Gaussian measure with
mean zero and variance $\sigma^{2}$. 
\end{theorem}

The assumptions of the theorem hold for the class of dynamical
systems discussed in \cite{young,CCS}.
For piece-wise expanding maps of the interval, a
stronger result was proved in \cite{chazottes-collet}.
Notice that this theorem immediately implies that almost-surely
$\mathcal{A}_{n}$ converges weakly to the Gaussian measure.

\bigskip

\begin{proof}
We first prove that
$$
\lim_{n\to\infty}\expectation\big(\ka{n}\big)=0\,.
$$
Let $B$ be a positive constant to be chosen large enough later on. We have
for any $g\in\van$ vanishing at $0$ and any $x$
$$
|g(x)|\le |x|\;.
$$
Therefore 
$$
\ka{n}
$$
\begin{equation}\label{mouton}
\le \sup_{g\in\van}\int_{-B}^{B}g \bigg( 
d\mathcal{A}_{n}-d\mathcal{N}\big(0,\sigma^{2}\big)\bigg)
+\int_{|y|>B}|y| d\mathcal{A}_{n}(y)
+\int_{|y|>B}|y| d\mathcal{N}\big(0,\sigma^{2}\big)(y)\;.
\end{equation}
We first estimate the expectation of the 
second term uniformly in $n$.
Since the correlations are absolutely summable, we get for any
$j$  
\begin{equation}\label{sume}
\expectation\big(S_{j}^{2}\big)^{1/2}
\le \Oun\sqrt j\;.
\end{equation}
Therefore, using Cauchy-Schwarz and Bienaym\'e-Chebychev inequalities we get
\begin{equation}\label{terme2}
\expectation\left(\int_{|y|>B}|y| d\mathcal{A}_{n}(y)\right)=
\frac{1}{D_{n}}\sum_{k=1}^{n}\frac{1}{k}\expectation\left(
\gchi_{[B,\infty[}\bigg(\big|S_{k}\big|/\sqrt
k\bigg)\frac{\big|S_{k}\big|}{\sqrt k}\right)\le \frac{\Oun}{B}\;\cdot
\end{equation}
In order to estimate the first term on the rhs of (\ref{mouton}), we observe that since $[-B,B]$ is
compact, we can apply Ascoli-Arzela theorem to conclude that for any
$\epsilon>0$ there is a number $r=r(\epsilon)$ and a finite sequence
$g_{1},\ldots,g_{r}$ of functions in $\van$ such that for any
$g\in\van$, there is at least one integer $1\le j\le r$ such that
$$
\sup_{|y|\le B}\big|g(y)-g_{j}(y)\big|\le\epsilon\;.
$$
Therefore
\begin{equation}\label{saute}
\sup_{g\in\van}\int_{-B}^{B}g \bigg( 
d\mathcal{A}_{n}-d\mathcal{N}\big(0,\sigma^{2}\big)\bigg)
\le \sup_{1\le j\le r(\epsilon)}\int_{-B}^{B}g_{j} \bigg( 
d\mathcal{A}_{n}-d\mathcal{N}\big(0,\sigma^{2}\big)\bigg)+2\epsilon\;.
\end{equation}
We now consider the $r$ sequences  of random variables 
$$
Y_{n,j}=\int_{-B}^{B}g_{j} \bigg( 
d\mathcal{A}_{n}-d\mathcal{N}\big(0,\sigma^{2}\big)\bigg)
$$
with $1\le j\le r$. 

We first estimate the variance of $Y_{n,j}$. 
Let the sequence of functions
$\big(K_{n,j}\big)$ of $n$ variables $x_{1},\ldots,x_{n}$ and $1\le j\le r$ be defined by
$$
K_{n,j}\big(x_{1},\ldots,x_{n}\big)=
\frac{1}{D_n} \sum_{k=1}^{n}\frac{1}{k}\left[
g_{j}\left(\frac{\sum_{l=1}^{k}u(x_{l})}{\sqrt k}\right)-\mathcal{N}(0,\sigma^{2})(g_j)
\right]
$$
where $\mathcal{N}(0,\sigma^{2})(\cdot)$ denotes the integration
against the Gaussian measure.
It is easy to verify that all these functions are separately Lipschitz
with respect to all their variables, and that the Lipschitz constant
with respect to the $q^{\mathrm{th}}$ variable is bounded by
$\Oun/(\sqrt{q}D_n)$ uniformly in $n$. Applying Devroye inequality (\ref{DEV}) we get
$$
\textup{var}\big(Y_{n,j})=
\textup{var}\big(K_{n,j}\big)
\le \frac{\Oun}{D_n^{2}}\sum_{q=1}^{n}\frac{1}{q} \le \frac{\Oun}{D_n}\;\cdot
$$ 
We now have using Cauchy-Schwarz inequality 
$$
\expectation\bigg(\sup_{1\le j\le r}Y_{n,j}\bigg)\le
\expectation\left(\sum_{j=1}^{r}\big|Y_{n,j}\big|\right)\le
\sum_{j=1}^{r}\expectation\left(\bigg|Y_{n,j}-\expectation\big(Y_{n,j}\big)
\bigg|\right)+\sum_{j=1}^{r}\bigg|\expectation\big(Y_{n,j}\big)\bigg|
$$
$$
\le\sum_{j=1}^{r}\textup{var}\big(Y_{n,j}\big)^{1/2}
+\sum_{j=1}^{r}\bigg|\expectation\big(Y_{n,j}\big)\bigg|
\leq 
\frac{r \Oun}{\sqrt{D_n}} + \sum_{j=1}^{r}\bigg|\expectation\big(Y_{n,j}\big)\bigg|
\;.
$$
By the central limit theorem (\ref{clt}), we have, for each  $j$ 
$$
\lim_{n\to\infty}\expectation\big(Y_{n,j}\big)=0\;,
$$
and therefore, from the above estimates, for a fixed $r(\epsilon)$ we have
$$
\limsup_{n\to\infty}\expectation\left(\sup_{1\le j\le
r(\epsilon)}\int_{-B}^{B}g_{j} \bigg(  
d\mathcal{A}_{n}-d\mathcal{N}\big(0,\sigma^{2}\big)\bigg)\right) \le 0\;.
$$
It now follows from (\ref{terme2}) and (\ref{saute}) that for any $\epsilon>0$ and any $B>0$
$$
0\le \limsup_{n\to\infty}\expectation\big(\ka{n}\big)\le
2\epsilon+\frac{\Oun}{B}+\int_{|y|>B}|y|\ d\mathcal{N}(0,\sigma^{2})(y)\;.
$$ 
Letting $B$ tend to infinity and $\epsilon$ to zero we get
$$
\lim_{n\to\infty}\expectation\big(\ka{n}\big)=0\,.
$$

We now estimate the variance of $\ka{n}$. 
Applying Devroye inequality (\ref{DEV}) as above to the function
$K_{n}$ of $n$ variables $x_{1},\ldots,x_{n}$ 
\begin{equation}
\label{bulle}
K_{n}\big(x_{1},\ldots,x_{n}\big)=
\end{equation}
$$
\sup_{g\in \van}\frac{1}{D_n} \sum_{j=1}^{n}\frac{1}{j}\left[
g\left(\frac{\sum_{l=1}^{j}u(x_{l})}{\sqrt
j}\right)-\mathcal{N}(0,\sigma^{2})(g)
\right]
$$
we get 
$$
\expectation\left(\bigg[\ka{n}-\expectation\big(\ka{n}\big)\bigg]^{2}\right)=
\textup{var}\big(K_{n}\big)
$$
$$
\le \frac{\Oun}{D_n^2}\sum_{j=1}^{n}\frac{1}{j} \le \frac{\Oun}{D_n}\;\cdot
$$ 
If for $0<\rho<1$ we define
$$
n_{k}=e^{k^{1+\rho}}
$$
we conclude that
$$
\sum_{k}
\expectation\left(\bigg[\ka{n_k}-\expectation\big(\ka{n_k}\big)\bigg]^{2}\right)
<\infty
$$
which implies by the B. L\'evi's theorem that 
$$
\lim_{k\to\infty}\big(\ka{n_k}-\expectation\big[\ka{n_k}\big]\big)=0
\quad\proba-\textup{almost surely}\,.
$$
We now observe that if $n_{k}<n\le n_{k+1}$ we have
$$
\left|\ka{n}-\ka{n_{k}}\right|
$$
$$
\le \frac{D_{n}-D_{n_k}}{D_{n}}\ka{n_{k}}+
\sup_{g\in \van}\frac{1}{D_n} \sum_{j=n_{k}+1}^{n}\frac{1}{j}\left[
g\left(\frac{S_{j}}{\sqrt
j}\right)-\mathcal{N}(0,\sigma^{2})(g)
\right]\;.
$$
The first term tends to zero almost surely by our previous estimates. We
now prove that the second term tends to zero almost surely. 
We have
$$
\sup_{g\in \mathscr L}\frac{1}{D_n} \sum_{j=n_{k}+1}^{n}\frac{1}{j}\left[
g\left(\frac{S_{j}}{\sqrt
j}\right)-\mathcal{N}(0,\sigma^{2})(g)
\right]
$$
$$
\le  \frac{1}{D_{n}}\sum_{j=n_{k}+1}^{n}\frac{1}{j}\left[
\frac{\big|S_{j}\big|}{\sqrt j}+\mathcal{N}(0,\sigma^{2})(|x|)
\right]
\le  \frac{1}{D_{n_{k}}}\sum_{j=n_{k}+1}^{n_{k+1}}\frac{1}{j}\left[
\frac{\big|S_{j}\big|}{\sqrt j}+\mathcal{N}(0,\sigma^{2})(|x|)
\right]\;.
$$
It follows easily from our choice of $(n_{k})$ that
$$
\lim_{k\to\infty}\frac{1}{D_{n_{k}}}\sum_{j=n_{k}+1}^{n_{k+1}}\frac{1}{j}
\mathcal{N}(0,\sigma^{2})(|x|)=0\;.
$$

We now prove the almost sure convergence to zero of the sequence
$$
T_k=\frac{1}{D_{n_{k}}}\sum_{j=n_{k}+1}^{n_{k+1}}\frac{\big|S_{j}\big|}{j^{3/2}}\,\cdot
$$
For this purpose we estimate the expectation of the square of $T_k$.
Using Cauchy-Schwarz inequality and (\ref{sume}) we obtain
$$
\expectation(T_k^2)\leq
\frac1{D_{n_k}^2}
\sum_{p,q= n_{k}+1}^{n_{k+1}} \frac{\expectation(S_p^2)^{1/2}}{p^{3/2}}
\frac{\expectation(S_q^2)^{1/2}}{q^{3/2}}
\leq
$$
$$
\frac{\left( \log n_{k+1}-\log n_k + \Oun\right)^2}{D_{n_{k}}^2}\leq \frac{\Oun}{k^2}\,\cdot
$$
It follows at once that $\expectation(T_k^2)$ is summable in $k$. The result now follows
using B. L\'evi's theorem. The theorem is proved.

\end{proof}

\noindent{\bf Remarks}. We note that the above proof also leads to an
estimate on the probability that $\ka{n}$ is larger than some given
number $\epsilon>0$. 

For a dynamical system $(\Omega,f)$ it often occurs that the invariant
measure is supported on an attractor which is a small subset of the phase
space $\Omega$.
When there exists a SRB measure, one would like to have Theorem \ref{t1}
almost-surely with respect to Lebesgue measure on $\Omega$.
Assuming that the stable foliation is absolutely continuous, and the
forward contraction is uniform and exponential along local stable manifolds (see
\cite{young} for several examples), it is 
sufficient to prove that 
$$
\lim_{n\to\infty} \left| 
K_n(x,f(x),\ldots,f^{n-1}(x))-K_n(\tilde{x},f(\tilde{x}),\ldots,f^{n-1}(\tilde{x}))
\right | =0
$$
where $K_n$ is defined by (\ref{bulle}) and $x,\tilde{x}$ belong to the same local stable manifold. This
follows at once from the definition of $K_n$ and the uniform
exponential contraction along local stable manifolds.

\appendix
\section{About a trigonometric series}\label{series}

For the convenience of the reader we prove in this appendix the following
(probably well-known) result on trigonometric series 
for which we have not been able to locate a reference.

\begin{lemma}\label{sin}
We have the following
$$
\sup_{m\in\mathbb{N},\omega\in[0,2\pi]}
\left| 
\sum_{k=1}^{m}\frac{\sin k\omega}{k}
\right|
<\infty\,.
$$
\end{lemma}

\begin{proof}
First observe that it is enough to assume that $\omega\in[0,\pi]$.
Now we have
$$
\sum_{k=1}^{m}\frac{\sin k\omega}{k}=\frac1{2} \int_0^\omega e^{is}\ \frac{1-e^{ims}}{1-e^{is}}\ ds
+ c.c.
$$
By an easy estimate, one gets
$$
\sup_{m\in\mathbb{N},\omega\in[0,\pi]}
\left| 
\int_{0}^\omega e^{is}\ \frac{1-e^{ims}}{1-e^{is}}\ ds-
\int_{0}^\omega e^{is}\ \frac{1-e^{ims}}{is}\ ds
\right|
<\infty\ .
$$
Finally, 
$$
\int_{0}^\omega e^{is}\ \frac{1-e^{ims}}{2is}\ ds + \, c.c.
$$
$$
=\int_0^\omega \frac{\sin s}{s}\ ds \ -
\int_0^\omega \frac{\sin(m+1) s}{s}\ ds 
=
-\int_\omega^{(m+1)\omega} \frac{\sin s}{s}\ ds\ .
$$
It is well-known that the modulus of this quantity
is uniformly bounded in $\omega$ and $m$.
\end{proof}

\section{On the uniform decay of correlations}\label{UDC}

In this appendix we prove a general result on decay of correlations
which may be useful in other contexts.
Consider a dynamical system on a phase space $\Omega$ given by a
measurable map $f$ from $\Omega$ to itself. Let $\mu$ be an ergodic
invariant measure. The 
decay of correlations is often proved in the following form:  
There is a non increasing sequence $(\gamma_{n})$ and two Banach
spaces $\bu$ and 
$\bd$ of measurable functions on $\Omega$ such that for any
functions $\psi_1\in\bu$ and $\psi_2\in\bd$, there is a constant
$C_{\psi_1,\psi_2}$ such that for any  integer $n$
\begin{equation}\label{gamman}
\left|\int \psi_1\circ f^{n}  \psi_2  \;d\mu -  
\int \psi_1  \;d\mu \int\psi_2  \;d\mu \right|
\le C_{\psi_1,\psi_2}\gamma_{n}\;.
\end{equation}
It is often useful to have  some information on the constant
$C_{\psi_1,\psi_2}$, in particular if  it can be bounded by a product of
norms of the two functions (and a uniform constant).
 It turns out that this apparently stronger result follows from
the previous estimate under the following natural assumptions.

\begin{itemize}
\item[\textup{{\em i)}}] The constant functions belong to $\bu$.
\item[\textup{{\em ii)}}] The integration with respect to $\mu$
defines a continuous linear functional on $\bu$.
\item[\textup{{\em iii)}}] The Koopman operator $U$ (of composition with $f$)
is continuous in $\bu$.
\item[\textup{{\em iv)}}] $\bd$ is contained in the dual of $\bu$
(duality with respect to the integration by $\mu$) with
a topology at least as fine as the dual norm topology. 
\end{itemize}

As will become clear  from the proof, the result below is due to the
special form of the correlation integral.  

\begin{theorem}\label{collet}
Assume the above properties i-iv),
and inequality  (\ref{gamman}) hold. Then 
there exists a constant $K$ such that for any integer $n$
and any $\psi_1\in \bu$, $\psi_2\in\bd$, we have
\begin{equation}\label{coco}
\left|\int \psi_1\circ f^{n}  \psi_2  \;d\mu -  
\int \psi_1  \;d\mu \int\psi_2  \;d\mu \right|
\le K\|\psi_1\|_{\bu}\|\psi_2\|_{\bd}\gamma_{n}\;.
\end{equation}
\end{theorem}

A frequent example is $\bu=L^{\infty}(\Omega,\;d\mu)$ while $\bd$ is a
space of more regular functions (functions of bounded variation,
Lipschitz or H\"older functions). 
In \cite{young}, $\bu=\bd$ is the space of H\"older continuous functions.
We give below a proof based on the
principle of uniform boundedness.  

\bigskip

\begin{proof}
We first deal with the easy case where for some integer $n_0$ we
have $\gamma_{n_0}=0$. For any $n>n_0$, using the identity
$\psi_1\circ f^n=(\psi_1\circ f^{n-n_0})\circ f^{n_0}$ 
and {\em iii)}, we conclude
that the correlation integral (left hand side of (\ref{gamman})) is equal to zero
for any $\psi_1\in\bu$ and $\psi_2\in\bd$. On the other hand, it follows
from {\em iv)} that there is a positive number $K_0$ such that
$$
\sup_{\|\psi_1\|_{\bu}\le 1\;,\;\|\psi_2\|_{\bd}\le 1}
\left|\int \psi_1 \;\psi_2 \;d\mu \right|\le K_0
$$
and (\ref{coco}) follows immediately  with
$$
K=K_0\;\sup_{n\;,\;\gamma_n>0}\|U^n\|_{\bu}\gamma_n^{-1}\;.
$$

We now assume $\gamma_n>0$ for any integer $n$. We first 
control the dependence on $\psi_1$ and for this
purpose we first fix $\psi_2\in\bd$. We then define a sequence of non
negative continuous  functions $(p^{\psi_2}_{n})$ on $\bu$ by
$$
p^{\psi_2}_{n}(\psi_1)=\gamma_{n}^{-1}\left|\int \psi_1\circ f^{n} 
 \psi_2  \;d\mu -  
\int \psi_1  \;d\mu \int\psi_2  \;d\mu \right|\;.
$$
We have obviously  for any integer $n$ and any $\psi_1$, $\psi_{1}'$ and
$\psi_1''$ belonging to $\bu$
$$
p^{\psi_2}_{n}(\psi_{1}'+\psi_{1}'')\le 
p^{\psi_2}_{n}(\psi_{1}')+p^{\psi_2}_{n}(\psi_{1}'')\quad
\hbox{\rm and }\quad 
p^{\psi_2}_{n}(\psi_1)=p^{\psi_2}_{n}(-\psi_1)\;.
$$
It follows immediately from  (\ref{gamman}) that  for each $\psi_1\in\bu$ we have
$$
\sup_{n}p^{\psi_2}_{n}(\psi_1)\le C_{\psi_1,\psi_2}<\infty\;.
$$
Therefore, we can apply the principle of uniform boundedness
\cite[Theorem 1.29 section III page 136]{K} to conclude that there is a finite
constant $D_{\psi_2}$ such that 
$$
\sup_{n,\|\psi_1\|_{\bu}\le 1}p^{\psi_2}_{n}(\psi_1)\le D_{\psi_2}\;.
$$
In other words, for any integer $n$, 
for any $\psi_1\in\bu$ and any $\psi_2\in\bd$ we have
\begin{equation}\label{toto}
\left|\int \psi_1\circ f^{n} 
 \psi_2  \;d\mu -  
\int \psi_1  \;d\mu \int\psi_2  \;d\mu \right|
\le D_{\psi_2}\|\psi_1\|_{\bu}\gamma_{n}\;.
\end{equation}
We shall now control the dependence in $\psi_2$.
Let $\Lambda=\integer\times B_{1}$ where $B_{1}$ is the closed  unit ball
of $\bu$. We define a family $(q_{\lambda})_{\lambda\in\Lambda}$ of
continuous, non-negative functions of $\bd$ by
$$
q_{(n,\psi_1)}(\psi_2)=\gamma_{n}^{-1}
\left|\int \psi_1\circ f^{n}   \psi_2  \;d\mu -  
\int \psi_1  \;d\mu \int\psi_2  \;d\mu \right|\;.
$$ 
We have immediately for any $\lambda\in\Lambda$ and for any $\psi_{2}$,
$\psi_{2}'$ and $\psi_2''$ in $\bd$
$$
q_{\lambda}(\psi_{2}'+\psi_{2}'')\le q_{\lambda}(\psi_{2}')
+q_{\lambda}(\psi_{2}'')\quad\hbox{\rm and } \quad
q_{\lambda}(\psi_2)=q_{\lambda}(-\psi_2)\;.
$$
Moreover it follows from (\ref{toto}) that for any $\psi_2\in\bd$
$$
\sup_{\lambda\in\Lambda}q_{\lambda}(\psi_2)\le D_{\psi_2}<\infty\;.
$$
We can apply as above  the principle  of uniform boundedness  to conclude
that there is a finite constant $K$ such that
$$
\sup_{\lambda\in\Lambda,\|\psi_2\|_{\bd}\le1}q_{\lambda}(\psi_2)\le K
$$
which immediately implies (\ref{coco}).
\end{proof}

In the case where $\gamma_n$ in (\ref{gamman}) is summable and
assumptions \textit{(i)-(iv)} hold, Theorem \ref{collet} implies Hypothesis
\ref{nuage} with $\bu=\bd$ being the space of $\eta$-H\"older continuous
functions ($0<\eta\leq 1$).

\section{A property of the density of the invariant measure for a class of 1D maps}\label{bing}

The purpose of this section is to prove that property (\ref{moutarde})
in Theorem \ref{graine} is indeed valid for maps on the interval satisfying the axioms
of \cite{young}. In other words the density of the absolutely continuous
invariant measure belongs to a Besov space (see \cite{stein} for definitions). 
In particular, quadratic maps for a set of parameters of positive
Lebesgue measure \cite{young} are included.
We refer the reader to \cite{young} (and \cite{CCS}) for notations and properties
of such dynamical systems and their associated tower maps. 

Recall that the density $\Phi$ of the SRB measure $\mu$ reads \cite{lulu,young}
\begin{equation}\label{Phi}
\Phi(y)=\sum_{j\geq 1} \sum_{k=0}^{R_j-1} a_{kj}(y) \gchi_{f^{k}(\Lambda_j)}(y)
\end{equation}
where we set, for any $k<R_j$ and for any $y\in f^{k}(\Lambda_j)$
$$
a_{kj}(y)=\frac{\varphi(y_{kj})}{f'^{k}(y_{kj})}
$$
where $y_{kj}$ is the unique point in $\Lambda_j$ satisfying
$f^k(y_{kj})=y$, and $\varphi$ is the density of the $f^R$-invariant measure.
It is convenient to assume that $a_{kj}$ vanishes outside $f^{k}(\Lambda_j)$.

We will use repeatedly the following properties coming from \cite{young}:
\begin{itemize}
\item[(i)] There exists $\theta>0$ such that $\sum_j e^{\theta R_j}\ |\Lambda_j|<\infty$.

\item[(ii)] There are constants $C>0$ and $\alpha\in]0,1[$ such that for all $j$ and all $k<R_j$
and any $y,y'$ in $f^{k}(\Lambda_j)$ 
$$
\left|
\frac{a_{kj}(y)}{a_{kj}(y')}-1
\right| \leq 
C\ \alpha^{s(y,y')}\,.
$$
We recall that $s(y,y')$ is the separation time of the orbits of $y$ and $y'$, see \cite{young};

\item[(iii)] There exists a constant $C>1$ such that for all $j$ and all $k<R_j$
and any $y$ in $f^{k}(\Lambda_j)$
$$
C^{-1}\ |\Lambda_j| \leq a_{kj}(y) |f^{k}(\Lambda_j)| \leq C\ |\Lambda_j|\,;
$$
\item[(iv)] Let $B:=\Vert f'\Vert_\infty>1$. For all $j$ and all $k<R_j$
and any $y$ in $f^{k}(\Lambda_j)$
$$
a_{kj}(y) \leq C B^{R_j-k}\ \frac{|\Lambda_j|}{|\Lambda|}\,\cdot
$$

\end{itemize}

Property (i) follows from the exponential tail for Markovian return times.
Property (ii) follows from the distortion bound in \cite{young}.
Property (iii) follows from (ii) and the fact that $f^k|\Lambda_j$ is
a diffeomorphism and $\varphi$ is bounded. Finally, property (iv) follows
from (iii) and the fact that $f^{R_j-k}(f^k(\Lambda_j))=\Lambda$
and $f^{R_j-k}|f^k(\Lambda_j)$ is a diffeomorphism.

We will use the following lemma.

\begin{lemma}\label{brouillard}
There exists a constant $C>0$ such that for any measurable set $A\in\mathbb{R}$ we have
$$
\mu(A)\leq C\ m(A)^{\varrho}
$$
where $m$ is Lebesgue measure and $\varrho=\min\{\theta/\log B, 1\}>0$ ($\theta$ and $B$ are defined
in (i) and (iv), respectively).
\end{lemma}

\begin{proof}
We have using H\"older inequality with $p=\log B/(\log B-\min\{\theta,\log B\})$
and $q=p/(p-1)=\varrho^{-1}$
$$
\mu(A)=\sum_{j\geq 1} \sum_{k=0}^{R_j-1} \int\ dy\,
 a_{kj}(y) \gchi_{A}(y)\ \gchi_{f^{k}(\Lambda_j)}(y)\leq
$$
$$
m(A)^{1/q}\
\sum_{j\geq 1} \sum_{k=0}^{R_j-1}
\left(
\int\ dy\,
 a_{kj}^{p}(y) \ \gchi_{f^{k}(\Lambda_j)}(y)
\right)^{1/p}\,.
$$
Using (iii), (iv) and (i) this is bounded above by
$$
\Oun\ m(A)^{1/q}\ \sum_{j\geq 1} \sum_{k=0}^{R_j-1} \ |\Lambda_j| \ B^{(R_j-k)(p-1)/p}
\leq \Oun \  m(A)^{\varrho}\,.
$$
The lemma is proved.
\end{proof}

The main result of this section is the following theorem.

\begin{theorem}
For an interval map satisfying hypotheses of \cite{young}, for any positive 
$\tau<\min\left\{\log\alpha^{-1}/(2\log B), \frac{1}{4} (\min \{1,\theta/\log B\})^3\right\}$,
there exists $C>0$ such that
$$
\int |\Phi(y) - \Phi(y-\delta)| \ dy \leq C\ |\delta|^\tau 
$$
for any $\delta\in\mathbb{R}$. In other words,
$\Phi$ belongs to the Besov space $\Lambda^{1,\infty}_\tau$ (see \cite{stein}).
\end{theorem}

\begin{proof}
It is enough to consider $0<\delta<1$. We have 
\begin{equation}\label{manger}
\int |\Phi(y) - \Phi(y-\delta)| \ dy \leq
$$
$$
\sum_{j\geq 1} \sum_{k=0}^{R_j-1} \int dy\, 
\left|
a_{kj}(y) \gchi_{f^{k}(\Lambda_j)}(y)
-
a_{kj}(y-\delta) \gchi_{f^{k}(\Lambda_j)}(y-\delta)
\right|\,.
\end{equation}
For a fixed $\delta>0$, we split the sum over $j$ and $k$ in (\ref{manger}) according to the condition
$\delta > |f^{k}(\Lambda_j)|/2$ and the complementary condition. The first sum is bounded above by
\begin{equation}\label{eueueu}
2\ \sum_{j\geq 1} \sum_{{\scriptstyle k=0} \atop
  |f^{k}(\Lambda_j)|<2\delta}^{R_j-1} \sup_{y\in f^{k}(\Lambda_j)}
a_{kj}(y)\  
|f^{k}(\Lambda_j)|=
$$
$$
2\ \sum_{j\geq 1} \sum_{{\scriptstyle k=0} \atop |f^{k}(\Lambda_j)|<2\delta}^{R_j-1} 
\left(\sup_{y\in f^{k}(\Lambda_j)} a_{kj}(y)\ |f^{k}(\Lambda_j)|\right)^{1-\tau}
\left(\sup_{y\in f^{k}(\Lambda_j)} a_{kj}(y)\right)^{\tau}\ |f^{k}(\Lambda_j)|^{\tau}
$$
$$
\leq \Oun\ \delta^{\tau}\ \sum_{j\geq 1} \sum_{k=0}^{R_j-1}
|\Lambda_j|\ B^{\tau(R_j-k)}\leq\Oun\ \delta^{\tau} 
\end{equation}

where this last inequality follows from (i), (iii), (iv).

\bigskip

Now we turn to the second sum, namely the sum over the indices $j$, $k$ satisfying
$|f^{k}(\Lambda_j)|\geq 2\delta$. This sum is bounded above 
by
$$
\sum_{j\geq 1} \sum_{{\scriptstyle k=0} \atop |f^{k}(\Lambda_j)|\geq 2\delta}^{R_j-1} 
\int dy\, 
a_{kj}(y) \left|\gchi_{f^{k}(\Lambda_j)}(y)-\gchi_{f^{k}(\Lambda_j)}(y-\delta)\right|\, +
$$
$$
\sum_{j\geq 1} \sum_{{\scriptstyle k=0} \atop |f^{k}(\Lambda_j)|\geq 2\delta}^{R_j-1} 
\int dy\, 
\left| a_{kj}(y)-a_{kj}(y-\delta)\right|\ 
\gchi_{f^{k}(\Lambda_j)}(y-\delta)\,.
$$
Since $f^{k}(\Lambda_j)$ is an interval and $|f^{k}(\Lambda_j)|\geq 2\delta$, we have
\begin{equation}\label{casino}
\int dy\, 
a_{kj}(y) \left|\gchi_{f^{k}(\Lambda_j)}(y)-\gchi_{f^{k}(\Lambda_j)}(y-\delta)\right|\leq
\end{equation}
$$
\Oun \ \sup_{y\in f^{k}(\Lambda_j)} a_{kj}(y)\ 
\int dy\, \left|\gchi_{f^{k}(\Lambda_j)}(y)-\gchi_{f^{k}(\Lambda_j)}(y-\delta)\right|\leq
$$ 
$$
\Oun \ \sup_{y\in f^{k}(\Lambda_j)} a_{kj}(y)\ \delta\,.
$$
On the other hand, using (iii), the same integral is bounded above by
$2\sup_{y\in f^{k}(\Lambda_j)} a_{kj}(y)\ |f^{k}(\Lambda_j)|\leq \Oun\ |\Lambda_j|$.
Proceeding as in (\ref{eueueu}) we obtain for the sum over $j$ and $k$ the upper
bound $\Oun \ \delta^{\tau}$.
We now estimate for each $j$ and $k$ the integral
$$
\int dy\, 
\left| a_{kj}(y)-a_{kj}(y-\delta)\right|\ 
\gchi_{f^{k}(\Lambda_j)}(y-\delta)=
$$
\begin{equation}\label{charolais}
\int_{(f^{k}(\Lambda_j))^c} dy\, 
\left| a_{kj}(y)-a_{kj}(y-\delta)\right|\ 
\gchi_{f^{k}(\Lambda_j)}(y-\delta) \, +
\end{equation}
\begin{equation}\label{entrecote}
\int_{f^{k}(\Lambda_j)} dy\, 
\left| a_{kj}(y)-a_{kj}(y-\delta)\right|\ 
\gchi_{f^{k}(\Lambda_j)}(y-\delta)\,.
\end{equation}
It is easy verify that the integral (\ref{charolais}) can be bounded above
like the integral (\ref{casino}).
For the integral (\ref{entrecote}) we have the obvious upper bound 
\begin{equation}\label{obvious}
2\sup_{y\in f^{k}(\Lambda_j)} a_{kj}(y)\ |f^{k}(\Lambda_j)|\leq \Oun\ |\Lambda_j|\,.
\end{equation}
Using (ii) this integral is also bounded above by
$$
\Oun\ \int_{f^{k}(\Lambda_j)} dy\, a_{kj}(y)\ \alpha^{s(y,y-\delta)}\ 
\gchi_{f^{k}(\Lambda_j)}(y-\delta) 
$$
where $s(y,y-\delta)$ is the separation time of the orbits of $y$ and $y-\delta$.
In order to estimate this integral we introduce a partition of $f^{k}(\Lambda_j)$
into four subsets defined by
$$
\mathbf{B}^1_{kj}= \left\{y\in f^{k}(\Lambda_j)\ \Big| \
  s(y,y-\delta)> \frac{\tau\log\delta}{\log\alpha}\right\}
$$
$$
\mathbf{B}^2_{kj}
= \left\{y\in f^{k}(\Lambda_j)\cap (\mathbf{B}^1_{kj})^c\ \Big|
R(f^{s(y,y-\delta)}(y)) > \sigma\ \log\delta^{-1}
\right\}\,
$$
where $R(\cdot)$ is the Markovian return-time function defined in \cite{young}, and $\sigma:=1/4\log B$.
$$
\mathbf{B}^3_{kj}= \left\{y\in f^{k}(\Lambda_j)\cap (\mathbf{B}^1_{kj})^c\cap\mathbf{B}^2_{kj}
\ \Big| \ |f^{s(y,y-\delta)}(y)-f^{s(y,y-\delta)}(y-\delta)| < \sqrt{\delta}
\right\}
$$
$$
\mathbf{B}^4_{kj}= \left\{y\in f^{k}(\Lambda_j)\cap (\mathbf{B}^1_{kj})^c\cap\mathbf{B}^2_{kj}
\ \Big| \ |f^{s(y,y-\delta)}(y)-f^{s(y,y-\delta)}(y-\delta)| \geq \sqrt{\delta}
\right\}\,.
$$
We will estimate the contribution of these four sequences of sets separately.
We have obviously
$$
\Oun\ \int_{\mathbf{B}^1_{kj}}\ a_{kj}(y)\ dy\, \alpha^{s(y,y-\delta)}\leq
\delta^{\tau}\ |f^{k}(\Lambda_j)|\ \sup_{y\in f^{k}(\Lambda_j)} a_{kj}(y)
$$
and therefore we can bound the sum over $k$ and $j$ using (iii) and (i).

To estimate the contribution of $\mathbf{B}^2_{kj}$ we introduce the set
$$
\mathcal{C}:=
\bigcup_{\ell=0}^{\lfloor \frac{\tau\log\delta}{\log\alpha}\rfloor}
\left\{y\ : \ f^{\ell}(y)\in \{R>\sigma\ \log\delta^{-1}\}\right\}\,.
$$
From the invariance of the SRB measure $\mu$ and Lemma \ref{brouillard} we have
$$
\mu(\mathcal{C})\leq
 \left\lfloor \frac{\tau\log\delta}{\log\alpha}\right\rfloor\mu\{R>\sigma\ \log\delta^{-1}\}
\leq C \ \frac{\tau\log\delta}{\log\alpha}\ m(R>\sigma\ \log\delta^{-1})^{\varrho}
\,.
$$
Now observe that 
$$
\mathbf{B}^2_{kj}\subset\mathcal{C}\cap f^{k}(\Lambda_j)
$$
which implies using (i) and Chebychev inequality that
$$
\int_{\mathbf{B}^2_{kj}} a_{kj}(y) \ dy \leq 
\sup_{y\in f^{k}(\Lambda_j)} a_{kj}(y) \
\mu(\mathcal{C}) \leq 
$$
$$
\Oun 
\sup_{y\in f^{k}(\Lambda_j)} a_{kj}(y)\
\frac{\tau\log\delta}{\log\alpha}\ m(R>\sigma\
\log\delta^{-1})^{\varrho} \leq
$$
$$
\Oun \ \log(\delta^{-1})\ \sup_{y\in f^{k}(\Lambda_j)} a_{kj}(y)
\ \delta^{\varrho\sigma\theta}\,.
$$
Using (iv) and interpolating with the bound (\ref{obvious}) we get
$$
\sum_{j\geq 1} \sum_{k=0}^{R_j-1} 
\int_{\mathbf{B}^2_{kj}} a_{kj}(y) \ dy \leq 
\Oun \log(\delta^{-1})
\delta^{\varrho^2 \sigma\theta}
= \Oun \log(\delta^{-1}) \delta^{\varrho^{3}/4 }\,.
$$
We now treat the integral over the set $\mathbf{B}_{kj}^3$.
We define the sets 
$$
\mathcal{D}_0=
\left\{ y\ :\ y\in\Lambda,\; d\left(
y,\cup_{j,R_j <\sigma\log\delta^{-1}} \partial\Lambda_j\right)<\sqrt{\delta}\right\}
$$
and
$$
\mathcal{D}=\bigcup_{\ell=0}^{\lfloor\sigma\log\delta^{-1}\rfloor} f^{-\ell}(\mathcal{D}_0)\,.
$$
From the invariance of the SRB measure $\mu$ we get
$$
\mu(\mathcal{D})\leq \sigma\log\delta^{-1}\ 
\mu(\mathcal{D}_0)\,.
$$
We now estimate
$$
\mu(\mathcal{D}_0)=\sum_{j\geq 1} \sum_{k=0}^{R_j-1}
\int_{f^k(\Lambda_j)\cap \mathcal{D}_0} a_{kj}(y)\ dy  \,.
$$
As we have done several times above, each integral in these sums has
two bounds. 
From the definition of $\mathcal{D}_0$ we have
$$
\int_{f^k(\Lambda_j)\cap \mathcal{D}_0} a_{kj}(y)\ dy \leq
2\sqrt{\delta} \ \#\{j\ :\ R_j<\sigma\log\delta^{-1}\}\ \sup_{y\in f^{k}(\Lambda_j)} a_{kj}(y)\,.
$$
Since $f^{R_j}(\Lambda_j)=\Lambda$ for all $j\geq 1$ we have for any $j\geq 1$ that
$$
B^{R_j} |\Lambda_j| \geq |\Lambda|\,.
$$ 
Therefore since the $\Lambda_j$'s are disjoint and their union is
$\Lambda$ we obtain for any integer $q\geq 1$
$$
\#\{j\ :\ R_j \leq q\}\leq B^q\,.
$$
It follows using (iv) that 
$$
\int_{f^k(\Lambda_j)\cap \mathcal{D}_0} a_{kj}(y)\ dy \leq
\Oun \ \delta^{1/4} B^{R_j-k} |\Lambda_j|\,.
$$
Therefore interpolating with the trivial bound (\ref{obvious}) as before one gets
$$
\mu(\mathcal{D}_0)\leq \Oun \ \delta^{\theta/(4\log B)}\,.
$$
We now observe that
$$
\mathbf{B}^3_{kj}\subset\mathcal{D}\cap f^{k}(\Lambda_j)
$$
which implies using (ii) that
$$
\int_{\mathbf{B}^3_{kj}} a_{kj}(y) \ dy \leq 
\sup_{y\in f^{k}(\Lambda_j)} a_{kj}(y) \
\mu(\mathcal{D}) \leq 
$$
$$
\Oun 
\sup_{y\in f^{k}(\Lambda_j)} a_{kj}(y)\
\log(\delta^{-1})\ \delta^{\theta/4\log B}\,.
$$
Using (i) and (iv) and interpolating with the bound (\ref{obvious}) we get
$$
\sum_{j\geq 1} \sum_{k=0}^{R_j-1} 
\int_{\mathbf{B}^3_{kj}} a_{kj}(y) \ dy 
\leq
\Oun\
\left(\log(\delta^{-1})\right)^{\theta/\log B}\ \delta^{\theta^2
  /4(\log B)^2}\,.
$$
Finally if $y\in \mathbf{B}^4_{kj}$ we have using $B=\Vert f'\Vert_\infty$ 
$$
\sqrt{\delta} \leq | f^{s(y,y-\delta)}(y)-f^{s(y,y-\delta)}(y-\delta)|
\leq \delta\ B^{s(y,y-\delta)}\,.
$$
This immediately implies that
$$
s(y,y-\delta)  \geq -\frac{\log\delta}{2\log B}\,\cdot
$$
Using this bound, properties (iii) and (i), we obtain
$$
\sum_{j\geq 1} \sum_{k=0}^{R_j-1} \
\int_{\mathbf{B}^4_{kj}}\ a_{kj}(y)\ \alpha^{s(y,y-\delta)}\ dy\, \leq 
$$
$$
\Oun \
\delta^{\log\alpha^{-1}/(2\log B)}\ \sum_{j\geq 1} \sum_{k=0}^{R_j-1}
|\Lambda_j| \leq \Oun \ \delta^{\log\alpha^{-1}/(2\log B)}\,.
$$

This ends the proof of the theorem. 

\end{proof}


\end{document}